\documentclass[10pt]{article}
\usepackage[T1]{fontenc}
\usepackage{lmodern}
\usepackage[margin=26mm]{geometry}
\usepackage{amsmath,amssymb,amsthm,mathtools}
\usepackage{xcolor}
\usepackage{indentfirst}
\usepackage{setspace}
\usepackage{enumerate}
\usepackage{needspace}
\usepackage{hyperref}
\usepackage{tikz}

\hypersetup{colorlinks=true,linkcolor=blue,citecolor=blue,urlcolor=blue,
  pdftitle={Integral abelianizations of higher spin Torelli subgroups}}

\newtheorem{theorem}{Theorem}[section]
\newtheorem{proposition}[theorem]{Proposition}
\newtheorem{lemma}[theorem]{Lemma}
\newtheorem{corollary}[theorem]{Corollary}
\theoremstyle{definition}

\theoremstyle{remark}
\newtheorem{remark}[theorem]{Remark}
\numberwithin{equation}{section}

\newcommand{\Z}{\mathbb Z}
\newcommand{\Q}{\mathbb Q}
\newcommand{\F}{\mathbb F}
\newcommand{\I}{\mathcal I}
\newcommand{\K}{\mathcal K}
\newcommand{\Mod}{\operatorname{Mod}}
\newcommand{\Sp}{\operatorname{Sp}}
\newcommand{\Tor}{\operatorname{Tor}}
\newcommand{\Arf}{\operatorname{Arf}}
\newcommand{\Sym}{\operatorname{Sym}}

\newcommand{\rk}{\operatorname{rank}}
\newcommand{\ab}{\mathrm{ab}}

\title{Integral abelianizations of higher spin Torelli subgroups}
\author{Qirong Yang}
\date{}

\begin{document}
\maketitle

\begin{abstract}
We determine the integral abelianizations of the Torelli subgroups
preserving an $r$-spin structure for $r\geqslant2$ and genus
$g\geqslant5$. Set $n=r/\gcd(r,2)$. For a surface with one boundary
component, the free rank agrees with that of the Torelli group,
and the torsion subgroup is the direct sum of the usual
Birman--Craggs--Johnson $2$-torsion and
$\mathbb{Z}/n^2\mathbb{Z}$. For a closed surface, an $r$-spin
structure exists precisely when $n\mid g-1$, and the same
description holds with the additional cyclic factor replaced by
$\mathbb{Z}/\gcd(n^2,g-1)\mathbb{Z}$. The additional cyclic torsion
is determined by the Casson--Morita homomorphism.
\end{abstract}

\noindent\textbf{Keywords:} Higher spin structures; Torelli groups;
Johnson homomorphism; integral abelianization.

\section{Introduction}\label{sec:introduction}

Let $\Sigma_{g,b}$ be a compact connected oriented surface
of genus $g$ with $b\in\{0,1\}$ boundary components.
Write $\Mod_{g,b}$ for its mapping class group and
$\I_{g,b}$ for the Torelli group, the kernel of the action
on $H_1(\Sigma_{g,b};\Z)$. For an integer $r\geqslant2$, an $r$-spin structure on
$\Sigma_{g,b}$ is a class in
$H^1(UT\Sigma_{g,b};\Z/r\Z)$ taking value $1$ on the
positively oriented fiber of the unit tangent bundle
$UT\Sigma_{g,b}\to\Sigma_{g,b}$. For an $r$-spin structure $\phi$,
let $\Mod_{g,b}[\phi]$ denote its stabilizer and define r-spin Torelli groups
\[
 \Gamma_{g,b}(r,\phi):=\I_{g,b}\cap\Mod_{g,b}[\phi].
\]

These subgroups have finite index in the Torelli group and
contain the Johnson kernel. They are independent of $\phi$;
we write $\Gamma_{g,b}(r)$ and $\Gamma_g(r):=\Gamma_{g,0}(r)$. By Putman's theorem~\cite[Theorem~B]{Putman}, these subgroups
have the same rational abelianizations as the Torelli groups
for $g\geqslant3$.
We show that their integral abelianizations detect additional finite torsion,
represented by Torelli commutators, whose order depends on
the spin level.

Spin mapping class groups occur as monodromy groups in
algebraic and flat geometry; see~\cite{SalterNotes,CalderonSalter}.
Sipe~\cite{Sipe82,Sipe86} studied the mapping class group
action on roots of the canonical bundle. Randal-Williams
\cite{RandalWilliams12} computed stable abelianizations
of higher spin mapping class groups. For classical spin
structures, Sierra~\cite{Sierra26} determined the first
homology of the stabilizers on surfaces with one boundary
component. Bianchi~\cite{Bianchi26} obtained finite
presentations of the even spin mapping
class group and recovered its abelianization. 

Johnson~\cite{Johnson80,Johnson85} determined the integral
abelianization of the Torelli group. Its free quotient is
the image of the first Johnson homomorphism, and its
torsion subgroup is detected by the Birman--Craggs--Johnson
homomorphism. Putman~\cite[Theorem~B]{Putman} proved that
finite-index subgroups containing the Johnson kernel have
the same rational abelianization as the Torelli group.
Thus the integral calculation for r-spin Torelli groups
reduces to the finite torsion subgroup.

The action of the Torelli group on $r$-spin structures
is determined by the Chillingworth homomorphism modulo $r$.
On $\I_{g,1}$, the Chillingworth homomorphism agrees,
up to sign, with twice the contraction of the Johnson
homomorphism~\cite[Theorem~2]{Johnson80}.
Calderon and Salter~\cite[Lemmas~5.12--5.14]{CalderonSalter}
describe closed r-spin Torelli groups by the corresponding
congruence. In the case of one boundary component, these
finite-level groups contain the Chillingworth subgroup,
the kernel of the integral Chillingworth homomorphism.
For $g\geqslant6$, Kosuge~\cite[Theorem~A]{Kosuge}
computed the rational abelianization of that subgroup.
On a closed surface, the image of the integral
Chillingworth subgroup is $\Gamma_g(2g-2)$.
The results below determine the commutator quotient and
the full finite torsion subgroup at each spin level,
including these coincident cases.

For $r\geqslant2$, set
\begin{equation}\label{eq:effective-level}
 n:=\frac{r}{\gcd(r,2)}.
\end{equation}
A closed surface admits an $r$-spin structure precisely
when $r\mid2g-2$, equivalently $n\mid g-1$.
For a surface with boundary there is no restriction on $r$.
We determine the integral abelianizations of r-spin Torelli
groups as follows.

\begin{theorem}\label{thm:main}
Let $g\geqslant5$, $r\geqslant2$, and let $n$ be as in
\eqref{eq:effective-level}. Then
\begin{equation}\label{eq:boundary-answer}
 H_1(\Gamma_{g,1}(r);\Z)
 \cong\Z^{\binom{2g}{3}}
 \oplus(\Z/2\Z)^{2g^2+g+1}\oplus\Z/n^2\Z.
\end{equation}
If $n\mid g-1$ and $m:=\gcd(n^2,g-1)$, then
\begin{equation}\label{eq:closed-answer}
 H_1(\Gamma_g(r);\Z)
 \cong\Z^{\binom{2g}{3}-2g}
 \oplus(\Z/2\Z)^{2g^2+g}\oplus\Z/m\Z.
\end{equation}
\end{theorem}

\begin{remark}\label{rem:classical-spin}
For $r=2$, a spin structure is equivalent to a quadratic
refinement of the intersection form on
$H_1(\Sigma_{g,b};\F_2)$.
Since the Torelli group acts trivially on this homology
group, it preserves every classical spin structure.
Consequently,
\[
 \Gamma_{g,b}(2)=\I_{g,b}
 \qquad(b=0,1),
\]
independently of the Arf invariant.
In this case $n=1$, and
Theorem~\ref{thm:main} recovers Johnson's computation of
the integral abelianization of the Torelli group
\cite{Johnson85}, which holds for every $g\geqslant3$.
\end{remark}

The cyclic factors are determined by quotients of
commutator subgroups. Write $G'=[G,G]$ and
$\gamma_3G=[G,G']$.

\begin{theorem}\label{thm:derived}
Under the hypotheses of Theorem~\ref{thm:main},
\[
 \gamma_3\I_{g,1}\leqslant\Gamma_{g,1}(r)',\qquad
 \I_{g,1}'/\Gamma_{g,1}(r)'\cong\Z/n^2\Z.
\]

If $n\mid g-1$, then
\[
 \gamma_3\I_g\leqslant\Gamma_g(r)',\qquad
 \I_g'/\Gamma_g(r)'\cong\Z/\gcd(n^2,g-1)\Z.
\]
\end{theorem}

The following corollary gives an explicit generator for the
cyclic summand $\Z/n^2\Z$ in the abelianization of
$\Gamma_{g,1}(r)$.

\begin{corollary}\label{cor:intro-cyclic-generator}
Let $g\geqslant5$, $r\geqslant2$, and $n=r/\gcd(r,2)$.
Let $f=T_{c_+}T_{c_-}^{-1}$ and
$h=T_{d_+}T_{d_-}^{-1}$ be genus-one bounding-pair maps
on $\Sigma_{g,1}$ whose oriented common homology classes
$a$ and $b$ satisfy $a\cdot b=1$.
Then the class of $[f,h]$ in
$H_1(\Gamma_{g,1}(r);\Z)$ generates a direct summand
isomorphic to $\Z/n^2\Z$.
\end{corollary}

We outline the proof of the main results.
Write $\I:=\I_{g,1}$, $\Gamma:=\Gamma_{g,1}(r)$,
$K:=\K_{g,1}$, and $n:=r/\gcd(r,2)$.
Since $\Gamma'\leqslant\I'\leqslant K\leqslant\Gamma$,
the Johnson homomorphism induces an exact sequence
\[
0\longrightarrow K/\Gamma'
\longrightarrow\Gamma^{\mathrm{ab}}
\xrightarrow{\tau}L\longrightarrow0,
\qquad L:=\tau(\Gamma).
\]
Putman's theorem~\cite[Theorem~B]{Putman} identifies $L\otimes\Q$ with
$H_1(\I;\Q)$ and shows that $K/\Gamma'$ is the finite
torsion subgroup of $\Gamma^{\mathrm{ab}}$.
Since $L$ is free abelian, this sequence splits.
It therefore remains to determine $K/\Gamma'$, which
fits into the exact sequence
\[
0\longrightarrow\I'/\Gamma'
\longrightarrow K/\Gamma'
\longrightarrow K/\I'\longrightarrow0.
\]
Here $K/\I'$ is the elementary $2$-group computed by
Johnson~\cite{Johnson85}. We must calculate the kernel $\I'/\Gamma'$
and determine whether this second sequence splits.

We first prove $[\I,K]\leqslant\Gamma'$ using
bounding-pair maps and the lantern relation.
The identity $\gamma_3\I=[\I,K]$ of
Faes--Massuyeau--Sato~\cite[Theorem~A]{FMS}
then shows that $\I/\Gamma'$ has nilpotency class
at most $2$. Since $x^n\in\Gamma$ for every $x\in\I$,
the identity $[x^n,y^n]=[x,y]^{n^2}$ in this quotient
gives $n^2(\I'/\Gamma')=0$.
Using the integral description of $\I'/\gamma_3\I$,
we prove that an element of $\I'$ belongs to $\Gamma'$
precisely when its normalized Casson--Morita value
$\lambda:=d_{\mathrm M}/8$ is divisible by $n^2$.
Thus $\lambda$ induces an isomorphism
$\I'/\Gamma'\cong\Z/n^2\Z$.
Together with the quotient map to $K/\I'$, it gives
\[
K/\Gamma'\xrightarrow{\;\cong\;}
(K/\I')\oplus\Z/n^2\Z,
\qquad
k\Gamma'\longmapsto
\bigl(k\I',\,\lambda(k)\bmod n^2\bigr).
\]
This proves that the torsion sequence splits for
every $n$, including even $n$.

For the closed surface, assume $n\mid g-1$.
The disk-pushing subgroup lies in $\Gamma$, and
capping reduces the free rank by $2g$.
On the torsion subgroup, capping kills the boundary-twist class, reducing the additional
cyclic factor to $\Z/\gcd(n^2,g-1)\Z$.

Section~\ref{sec:preliminaries} gives notation and describes
the r-spin Torelli groups.
Sections~\ref{sec:free} and~\ref{sec:two-torsion} identify
the free quotient and the $2$-torsion $K/\I'$, respectively.
Section~\ref{sec:third-term} proves
$\gamma_3\I\leqslant\Gamma'$ by showing that separating
twists have central images in $\I/\Gamma'$.
Section~\ref{sec:commutator} computes the commutator quotient,
and Section~\ref{sec:splitting} splits the torsion sequence
and describes generators of the abelianization.
Section~\ref{sec:capping} deduces the closed case
by capping the boundary.
Section~\ref{sec:low-genus} gives obstructions in genera
$3$ and $4$.

\section{Preliminaries}\label{sec:preliminaries}

\subsection{Surfaces and Johnson homomorphisms}

Mapping classes and isotopies fix the boundary pointwise.
Write $\Sigma_g:=\Sigma_{g,0}$ and $\Mod_g:=\Mod_{g,0}$.
Let $H=H_1(\Sigma_{g,b};\Z)$, with intersection form
$\omega$. Capping identifies the homology groups for
$b=0$ and $b=1$. The Torelli group is
\[
 \I_{g,b}:=\ker\bigl(\Mod_{g,b}\longrightarrow\Sp(H,\omega)\bigr).
\]
We assume $g\geqslant3$ unless stated otherwise and use
$[x,y]=xyx^{-1}y^{-1}$. The lower central series is defined
by $\gamma_1G=G$ and $\gamma_{j+1}G=[G,\gamma_jG]$.
For an abelian group $A_0$, write $\Tor(A_0)$ for its
torsion subgroup. Tensor products are over $\Z$.

Fix a symplectic basis $a_1,b_1,\ldots,a_g,b_g$ of $H$,
with $\omega(a_i,b_j)=\delta_{ij}$, and put
\[
 \Omega:=\sum_i a_i\wedge b_i,\quad V:=\Lambda^3H,
 \quad\iota:H\longrightarrow V,\quad x\longmapsto\Omega\wedge x.
\]
The map $\iota$ is injective with primitive image. Indeed,
for each basis vector $x$, choose a handle $j$ different
from its handle. The coefficient of
$a_j\wedge b_j\wedge x$ in $\iota(v)$ recovers the
coefficient of $x$ in $v$. These coefficient maps give
a retraction of $\iota$. Consequently
\[
 A:=V/\iota(H)\quad\text{is free},\qquad
 \rk V=\binom{2g}{3},\quad\rk A=\binom{2g}{3}-2g.
\]
The first Johnson homomorphisms are surjective maps
\[
 \tau_{g,1}:\I_{g,1}\twoheadrightarrow V,
 \qquad\tau_g:\I_g\twoheadrightarrow A.
\]
Their kernels are denoted by $\K_{g,1}$ and $\K_g$.
Our normalization on a bounding pair $c,d$ cobounding
a subsurface $S$ is
\begin{equation}\label{eq:bp-johnson}
 \tau_{g,1}(T_cT_d^{-1})=[c]\wedge\Omega_S,
\end{equation}
where $\Omega_S$ is the symplectic bivector of its genus
part; see~\cite{Johnson80} and \cite[\S2.4]{Putman}.

Define
\begin{equation}\label{eq:contraction}
 C(x\wedge y\wedge z)
 :=\omega(x,y)z+\omega(y,z)x+\omega(z,x)y,
 \qquad W:=\ker C.
\end{equation}
The map $C:V\to H$ is surjective: wedging a basis vector
with a symplectic pair on another handle gives a preimage.
Moreover,
\begin{equation}\label{eq:contraction-iota}
 C\iota(x)=(g-1)x.
\end{equation}
For $n\geqslant1$, let
\[
C_n:V\longrightarrow H/nH,
\qquad
C_n(v)=Cv\bmod nH.
\]
If $n\mid g-1$, this map vanishes on $\iota(H)$
and induces a map $C_n:A\to H/nH$, which we denote
by the same symbol.

\subsection{$r-$spin Torelli groups}

Choose a Riemannian metric on $\Sigma_{g,b}$, and let
$p:UT\Sigma_{g,b}\to\Sigma_{g,b}$ be the projection
of its unit tangent bundle. An $r$-spin structure is
a class $\phi\in H^1(UT\Sigma_{g,b};\Z/r\Z)$ taking value
$1$ on the oriented fiber. This cohomological definition is equivalent to the
winding-number definition of Humphries and Johnson
\cite{HJ89}; see also \cite[Theorem~2.2]{Hamenstadt}.
The Gysin sequence shows that, if an $r$-spin structure
$\phi_0$ exists, every $r$-spin structure can be written
uniquely as $\phi_0+p^*\eta$, where
$\eta\in H^1(\Sigma_{g,b};\Z/r\Z)$.
For $b=1$, $r-$spin structures always exist; for $b=0$ they exist for $r\mid2g-2$.

\begin{proposition}\label{prop:spin-kernel}
Let $r\geqslant2$ and $n=r/\gcd(r,2)$. The subgroup
$\Gamma_{g,b}(r,\phi)$ is independent of $\phi$, and

\begin{equation}\label{eq:spin-kernels}
\Gamma_{g,1}(r)=\ker(C_n\circ\tau_{g,1}),
\qquad
\Gamma_g(r)=\ker(C_n\circ\tau_g),
\end{equation}
where the closed formula assumes $n\mid g-1$.
\end{proposition}
\begin{proof}
For $f\in\I_{g,b}$, write
\[
 f^*\phi-\phi=p^*\xi_\phi(f),\qquad
 \xi_\phi(f)\in H^1(\Sigma_{g,b};\Z/r\Z).
\]
The class is unique by the Gysin sequence. Replacing
$\phi$ by $\phi+p^*\eta$ does not change it, because
$f$ acts trivially on $H^1(\Sigma_{g,b};\Z/r\Z)$.
Thus $\xi_\phi$ is independent of $\phi$.
Moreover, $f$ preserves $\phi$ if and only if $\xi_\phi(f)=0$.

On the surface with one boundary component, the
Chillingworth formula identifies this change, using the
intersection pairing, with
\[
 \xi_\phi(f)=\pm2C\tau_{g,1}(f)\pmod r
\]
\cite[Theorem~2]{Johnson80}. 

For $v\in H$, the condition $2v\in rH$ is equivalent to
$v\in nH$. Thus the formula above gives
\[
 f\in\Gamma_{g,1}(r)
 \quad\Longleftrightarrow\quad
 C_n\tau_{g,1}(f)=0.
\]
For the closed surface, $r\mid 2g-2$ implies $n\mid g-1$.
Since $C\iota(x)=(g-1)x$, the reduction of $C$ modulo $n$
is well defined on $A=V/\iota(H)$.
Given $f\in\I_g$, choose a lift
$\widetilde f\in\I_{g,1}$ and restrict $\phi$ to
$\Sigma_{g,1}$. Applying the preceding formula gives
\[
 f\in\Gamma_g(r)
 \quad\Longleftrightarrow\quad
 C\tau_{g,1}(\widetilde f)\in nH
 \quad\Longleftrightarrow\quad
 C_n\tau_g(f)=0.
\]

Both kernel descriptions show that the r-spin Torelli
groups contain the Johnson kernels. The equivariance
of the Johnson homomorphisms and $C_n$ shows that
these groups are normal in the corresponding mapping
class groups. 
\end{proof}

\begin{remark}\label{rem:chillingworth}
On $\Sigma_{g,1}$, the integral Chillingworth homomorphism
is, up to sign, $2C\tau_{g,1}$.
Its kernel is therefore
\[
 \operatorname{Ch}_{g,1}
 :=\ker(C\tau_{g,1}),
\]
whereas $\Gamma_{g,1}(r)$ is defined by the congruence
$C\tau_{g,1}\equiv0\pmod n$.
The homomorphism $C\tau_{g,1}$ induces an isomorphism
\[
 \Gamma_{g,1}(r)/\operatorname{Ch}_{g,1}
 \cong nH.
\]
Thus $\operatorname{Ch}_{g,1}$ has infinite index in
each finite-level spin Torelli group, and
\[
 \operatorname{Ch}_{g,1}
 =\bigcap_{s\geqslant2}\Gamma_{g,1}(s).
\]
For a closed surface, the Chillingworth subgroup
$\operatorname{Ch}_g$ is defined as the capping image of
$\operatorname{Ch}_{g,1}$~\cite[Definition~3.2]{Kosuge}.
It has finite index in $\I_g$ and satisfies
\[
 \operatorname{Ch}_g=\Gamma_g(2g-2).
\]
It is contained in every admissible $\Gamma_g(r)$;
when $g-1$ is odd, it also equals $\Gamma_g(g-1)$.
These closed-surface identifications follow from
$C\iota=(g-1)\operatorname{id}$ and the Johnson
homomorphism on the disk-pushing subgroup, as explained
in Section~\ref{sec:capping}.
\end{remark}

\subsection{The Casson--Morita homomorphism}

Let $d_{\mathrm M}:\K_{g,1}\to\Z$ be Morita's core of
the Casson invariant~\cite{Morita89,Morita91}, normalized by
\[
 d_{\mathrm M}(T_c)=4h(h-1)
\]
when $c$ bounds genus $h$. Separating twists generate
$\K_{g,1}$~\cite[Theorem~A]{Putman}, and $h(h-1)$ is even.
Hence
\begin{equation}\label{eq:lambda-definition}
 \lambda:=d_{\mathrm M}/8:\K_{g,1}\longrightarrow\Z
\end{equation}
is integral. Morita's commutator formula is
\begin{equation}\label{eq:morita-commutator}
 \lambda([f,h])
 =\omega(C\tau_{g,1}(f),C\tau_{g,1}(h))
 \qquad(f,h\in\I_{g,1});
\end{equation}
see~\cite[Proposition~3.2]{FMS} for these conventions.
In particular, $\lambda(\I_{g,1}')=\Z$: choose $f,h$ with
$C\tau_{g,1}(f)=a_1$ and $C\tau_{g,1}(h)=b_1$.

\section{The free part}\label{sec:free}

For calculations on $\Sigma_{g,1}$, write
\[
 \Sigma:=\Sigma_{g,1},\qquad
 \I:=\I_{g,1},\qquad K:=\K_{g,1},\qquad
 \Gamma:=\Gamma_{g,1}(r),\qquad\tau:=\tau_{g,1},
 \qquad\alpha:=C\tau.
\]
Put $L:=\tau(\Gamma)$. Since $\tau$ is surjective,
Proposition~\ref{prop:spin-kernel} gives
\[
 L=\ker(C_n:V\to H/nH)=C^{-1}(nH),
 \qquad\Gamma'\leqslant\I'\leqslant K\leqslant\Gamma.
\]
The map $C_n\tau$ has finite target, so $\Gamma$ has
finite index in $\I$. The same argument applies to
$\Gamma_g(r)$ when $n\mid g-1$.
Since $\Gamma$ is normal in $\I$, so is $\Gamma'$.
The Johnson homomorphism therefore induces
\begin{equation}\label{eq:free-sequence}
 0\longrightarrow K/\Gamma'\longrightarrow\Gamma^{\ab}
 \xrightarrow{\tau}L\longrightarrow0.
\end{equation}
We use the following form of Putman's theorem, specialized
to the Torelli group of the whole surface.

\begin{theorem}[Putman, {\cite[Theorem~B]{Putman}}]\label{thm:putman}
Let $g\geqslant3$ and $b\in\{0,1\}$. If
$F\leqslant\I_{g,b}$ has finite index and contains
$\K_{g,b}$, the Johnson homomorphism induces an isomorphism
\[
 H_1(F;\Q)\xrightarrow{\cong}
 \tau_{g,b}(\I_{g,b})\otimes\Q,
\]
where $\tau_{g,0}:=\tau_g$.
\end{theorem}

\begin{corollary}\label{cor:free-quotient}
For $g\geqslant3$, the Johnson homomorphism and inclusion
induce isomorphisms
\[
 H_1(\Gamma;\Q)\cong L\otimes\Q
 \cong H_1(\I;\Q)\cong V\otimes\Q.
\]
The free quotient $L$ has rank $\binom{2g}{3}$.
If $n\mid g-1$, the closed Johnson image
$\tau_g(\Gamma_g(r))=\ker(C_n:A\to H/nH)$
is free of rank $\binom{2g}{3}-2g$.
\end{corollary}
\begin{proof}
The r-spin Torelli groups have finite index and contain
the Johnson kernels, so Theorem~\ref{thm:putman} applies.
The subgroup $L\subseteq V$ has finite index because
$V/L\cong H/nH$. Hence $L\otimes\Q=V\otimes\Q$,
and $L$ has the same rank as $V$.
In the closed case, $\ker C_n\subseteq A$ also has finite
index and therefore has the same rank as $A$.
\end{proof}

\begin{lemma}\label{lem:finite-torsion}
Let $g\geqslant3$, $b\in\{0,1\}$, and let
$F\leqslant\I_{g,b}$ have finite index and contain
$\K_{g,b}$. Then $\K_{g,b}/F'$ is the finite torsion
subgroup of $F^{\ab}$.
\end{lemma}
\begin{proof}
Since $F'\leqslant\K_{g,b}\leqslant F$, the Johnson
homomorphism gives an exact sequence
\[
 0\longrightarrow\K_{g,b}/F'
 \longrightarrow F^{\ab}
 \longrightarrow\tau_{g,b}(F)\longrightarrow0.
\]
The group $\tau_{g,b}(F)$ is free abelian, so every
torsion element of $F^{\ab}$ belongs to the kernel.
By Theorem~\ref{thm:putman}, the last map becomes an
isomorphism over $\Q$; thus every element of its kernel
is torsion. Finally, $\I_{g,b}$ is finitely generated
\cite{Johnson83}, and $F$ has finite index in it.
Therefore $F^{\ab}$ is finitely generated and its torsion
subgroup $\K_{g,b}/F'$ is finite.
\end{proof}

\section{The Birman--Craggs--Johnson quotient}
\label{sec:two-torsion}

Let $\overline H=H/2H$, and let $\mathcal Q$ be the set
of quadratic refinements of its intersection form. Thus
$q:\overline H\to\F_2$ belongs to $\mathcal Q$ if
\[
 q(x+y)=q(x)+q(y)+\omega(x,y).
\]
For $x\in\overline H$, write $\overline x(q)=q(x)$.
Let $B_{\leqslant k}$ be the space of Boolean polynomials
of degree at most $k$ in $\overline{a_i},\overline{b_i}$,
viewed as functions on $\mathcal Q$. Put
\[
 \Arf:=\sum_i\overline{a_i}\,\overline{b_i},\qquad
 \mathcal Q_0:=\{q\in\mathcal Q:\Arf(q)=0\}.
\]
Write $B_{\leqslant2}^0$ for the image of $B_{\leqslant2}$
under restriction to $\mathcal Q_0$. Johnson's theorem
\cite{Johnson85} gives
\begin{equation}\label{eq:boolean-dimensions}
 \begin{gathered}
 \dim_{\F_2}B_{\leqslant2}=2g^2+g+1,\qquad
 B_{\leqslant2}^0=B_{\leqslant2}/\langle\Arf\rangle,\\
 \dim_{\F_2}B_{\leqslant2}^0=2g^2+g.
 \end{gathered}
\end{equation}
The restriction of the Birman--Craggs--Johnson
homomorphism $\beta:\I\to B_{\leqslant3}$ to $K$
induces isomorphisms
\begin{equation}\label{eq:johnson-torsion}
 K/\I'\xrightarrow[\beta]{\cong}B_{\leqslant2},\qquad
 \K_g/\I_g'\cong B_{\leqslant2}^0.
\end{equation}
If $c$ bounds genus $h$ with symplectic basis
$u_1,v_1,\ldots,u_h,v_h$, then
\begin{equation}\label{eq:bcj-twist}
 \beta(T_c)(q)=\sum_{i=1}^h q(u_i)q(v_i).
\end{equation}
In particular, the boundary twist has BCJ image $\Arf$.
See also~\cite[\S2.1]{FMS}.

The torsion subgroup in \eqref{eq:free-sequence} fits into
\begin{equation}\label{eq:torsion-extension}
 0\longrightarrow\I'/\Gamma'\longrightarrow K/\Gamma'
 \longrightarrow K/\I'\longrightarrow0.
\end{equation}
Equation~\eqref{eq:johnson-torsion} identifies its quotient
with the elementary $2$-group computed by Johnson.
The kernel is calculated in Section~\ref{sec:commutator}.
Section~\ref{sec:splitting} then determines the extension,
including the case in which $n$ is even.
\section{Centrality of separating twists}
\label{sec:third-term}

Set $N:=\I/\Gamma'$. We prove that every separating
twist has central image in $N$. We first determine the
contracted Johnson images of Torelli elements supported
on a subsurface.

\begin{lemma}\label{lem:local-alpha}
Let $P\subseteq\Sigma$ be a connected subsurface of genus
at least $2$ with one boundary component. If $\I(P)$ denotes
its Torelli group extended by the identity, then
\[
 \alpha(\I(P))=H_1(P;\Z)\subseteq H.
\]
\end{lemma}
\begin{proof}
Naturality of the Johnson homomorphism gives
$\alpha(\I(P))\subseteq H_1(P;\Z)$.
Since $P$ has genus at least $2$, it contains a genus-one
bounding pair $c,d$ with primitive common class $a$.
Choose a symplectic pair $u,v$ in the genus-one subsurface
between $c$ and $d$. Then $a$ is orthogonal to $u,v$, and
\eqref{eq:bp-johnson} gives
\[
 \tau(T_cT_d^{-1})=a\wedge u\wedge v,
 \qquad\alpha(T_cT_d^{-1})=a.
\]
Mapping classes supported on $P$ realize $\Sp(H_1(P;\Z))$,
which is transitive on primitive vectors.
Conjugating this bounding-pair map therefore shows that
every primitive vector belongs to $\alpha(\I(P))$.
These vectors generate $H_1(P;\Z)$, proving the reverse
inclusion.
\end{proof}

We shall repeatedly use the following supported correction.
If $f,u\in\I$, $\alpha(u)=\alpha(f)$, and $u$ commutes
with a separating twist $T_c$, then
$fu^{-1}\in\ker\alpha\leqslant\Gamma$ and $T_c\in K\leqslant\Gamma$.
Consequently,
\begin{equation}\label{eq:correction}
 [f,T_c]=[fu^{-1},T_c]\in\Gamma'.
\end{equation}
The equality follows directly from $u^{-1}T_cu=T_c$.

\begin{proposition}\label{prop:centrality}\label{prop:third-term}
For $g\geqslant5$ and every $n\geqslant1$,
\begin{equation}\label{eq:centrality}
 [\I,K]\leqslant\Gamma'.
\end{equation}
In particular, $\gamma_3\I\leqslant\Gamma'$ and $N$ has
nilpotency class at most $2$.
\end{proposition}
\begin{proof}
Separating twists generate $K$ \cite[Theorem~A]{Putman},
so it suffices to show that each has central image in $N$.
A boundary-parallel twist is already central in
$\Mod_{g,1}$. We may therefore assume that $c$ is an
essential separating curve not parallel to the boundary.

\emph{Both sides have genus at least $2$.}
Fix $f\in\I$. The two components of $\Sigma\setminus c$
give symplectic homology summands, and we write
\[
 \alpha(f)=v_{\mathrm L}+v_{\mathrm R}
\]
according to this decomposition.
On the side containing $\partial\Sigma$, cut along an arc
joining $c$ to $\partial\Sigma$ and take a smaller subsurface
of the same genus with one boundary component. Its first
homology is the same summand of $H$.
Lemma~\ref{lem:local-alpha} gives $u_{\mathrm L},u_{\mathrm R}$
supported on the two sides with images $v_{\mathrm L},v_{\mathrm R}$.
Their product $u$ satisfies $\alpha(u)=\alpha(f)$ and
commutes with $T_c$. Equation~\eqref{eq:correction} gives
$[f,T_c]\in\Gamma'$.

\emph{The side away from $\partial\Sigma$ has genus $1$.}
Write $c=\delta_1=\partial P_1$, where $P_1$ is a
one-holed torus. In its complement choose disjoint
one-holed tori $P_2,P_3$ so that $P_1,P_2,P_3$ lie in
a genus-three subsurface $P_0$ with one boundary component.
Put $\delta_i=\partial P_i$ for $0\leqslant i\leqslant3$.
The four-holed sphere between these tori has lantern relation
\begin{equation}\label{eq:lantern}
 T_{\delta_0}T_{\delta_1}T_{\delta_2}T_{\delta_3}
 =T_{\eta_{12}}T_{\eta_{13}}T_{\eta_{23}},
\end{equation}
where $\eta_{ij}$ bounds the genus-two subsurface containing
$P_i$ and $P_j$; see Figure~\ref{fig:lantern} and
\cite[\S3]{ChenChenLanier}.
The genera on the two sides of $\eta_{ij}$ are $2$ and
$g-2$, while those on the two sides of $\delta_0$ are
$3$ and $g-3$. Thus all these twists have central images
in $N$ by the first case.

\begin{figure}[htbp]
\centering
\begin{tikzpicture}[x=.0125cm,y=-.0125cm,
  line cap=round,line join=round,every node/.style={font=\small}]
\definecolor{lanternred}{HTML}{C42F3D}
\definecolor{lanternblue}{HTML}{245FC7}
\definecolor{lanternpurple}{HTML}{96359E}
\draw[black,line width=.85pt] (330,300) -- (330,245) .. controls (330,105) and (550,105) .. (550,245) -- (550,300);
\draw[black,line width=.85pt] (250,355) -- (195,355) .. controls (55,355) and (55,575) .. (195,575) -- (250,575);
\draw[black,line width=.85pt] (330,630) -- (330,690) .. controls (330,830) and (550,830) .. (550,690) -- (550,630);
\draw[black,line width=.85pt] (630,355) -- (925,355) .. controls (1065,355) and (1065,575) .. (925,575) -- (630,575);
\draw[black,line width=.85pt] (330,300) .. controls (330,335) and (300,355) .. (250,355);
\draw[black,line width=.85pt] (550,300) .. controls (550,335) and (580,355) .. (630,355);
\draw[black,line width=.85pt] (250,575) .. controls (300,575) and (330,600) .. (330,630);
\draw[black,line width=.85pt] (630,575) .. controls (580,575) and (550,600) .. (550,630);
\draw[black,line width=.85pt,densely dotted] (330,300) .. controls (385,280) and (495,280) .. (550,300);
\draw[black,line width=.85pt] (330,300) .. controls (385,320) and (495,320) .. (550,300);
\draw[black,line width=.85pt,densely dotted] (250,355) .. controls (230,410) and (230,520) .. (250,575);
\draw[black,line width=.85pt] (250,355) .. controls (270,410) and (270,520) .. (250,575);
\draw[black,line width=.85pt,densely dotted] (330,630) .. controls (385,610) and (495,610) .. (550,630);
\draw[black,line width=.85pt] (330,630) .. controls (385,650) and (495,650) .. (550,630);
\draw[black,line width=.85pt,densely dotted] (630,355) .. controls (610,410) and (610,520) .. (630,575);
\draw[black,line width=.85pt] (630,355) .. controls (650,410) and (650,520) .. (630,575);
\draw[lanternred,line width=.75pt,densely dotted] (308.75,340.625) .. controls (385,430) and (500,510) .. (571.25,591.25);
\draw[lanternblue,line width=.75pt,densely dotted] (308.75,591.25) .. controls (385,500) and (495,425) .. (571.25,340.625);
\draw[lanternpurple,line width=.75pt,densely dotted] (308.75,340.625) .. controls (382,399) and (503,399) .. (571.25,340.625);
\draw[lanternpurple,line width=.75pt,densely dotted] (571.25,591.25) .. controls (495,542) and (385,542) .. (308.75,591.25);
\draw[lanternred,line width=.75pt] (308.75,340.625) .. controls (407,390) and (510,501) .. (571.25,591.25);
\draw[lanternblue,line width=.75pt] (308.75,591.25) .. controls (415,534) and (526,420) .. (571.25,340.625);
\draw[lanternpurple,line width=.75pt] (571.25,340.625) .. controls (608,414) and (608,516) .. (571.25,591.25);
\draw[lanternpurple,line width=.75pt] (308.75,591.25) .. controls (346,517) and (346,411) .. (308.75,340.625);
\draw[black,line width=.85pt] (415,211) .. controls (415,244) and (465,244) .. (465,211);
\draw[black,line width=.85pt] (419,224) .. controls (428,205) and (452,205) .. (461,224);
\draw[black,line width=.85pt] (140,456) .. controls (140,489) and (190,489) .. (190,456);
\draw[black,line width=.85pt] (144,469) .. controls (153,450) and (177,450) .. (186,469);
\draw[black,line width=.85pt] (415,706) .. controls (415,739) and (465,739) .. (465,706);
\draw[black,line width=.85pt] (419,719) .. controls (428,700) and (452,700) .. (461,719);
\draw[black,line width=.85pt] (722,456) .. controls (722,489) and (772,489) .. (772,456);
\draw[black,line width=.85pt] (726,469) .. controls (735,450) and (759,450) .. (768,469);
\draw[black,line width=.85pt] (885,456) .. controls (885,489) and (935,489) .. (935,456);
\draw[black,line width=.85pt] (889,469) .. controls (898,450) and (922,450) .. (931,469);
\draw[black,line width=.85pt] (1012,465) .. controls (1012,473.837) and (1004.84,481) .. (996,481) .. controls (987.163,481) and (980,473.837) .. (980,465) .. controls (980,456.163) and (987.163,449) .. (996,449) .. controls (1004.84,449) and (1012,456.163) .. (1012,465);
\node[text=black] at (582,211) {$P_1$};
\node[text=black] at (105,334) {$P_2$};
\node[text=black] at (583,729) {$P_3$};
\node[text=black] at (576,300) {$\delta_1$};
\node[text=black] at (246,607) {$\delta_2$};
\node[text=black] at (576,633) {$\delta_3$};
\node[text=black] at (657,340) {$\delta_0$};
\node[text=lanternred] at (539,513) {$\eta_{23}$};
\node[text=lanternblue] at (509,401) {$\eta_{12}$};
\node[text=lanternpurple] at (365,463) {$\eta_{13}$};
\node[text=black] at (828,465) {$\cdots$};
\node[text=black] at (847,325) {genus $g-3\geqslant2$};
\node[text=black] at (996,513) {$\partial\Sigma$};
\end{tikzpicture}
\caption{The lantern used in the proof. 
The right component has genus $g-3\geqslant2$ and contains
$\partial\Sigma$. The curve $\eta_{ij}$ bounds a genus-two
subsurface containing $P_i$ and $P_j$.
}
\label{fig:lantern}
\end{figure}
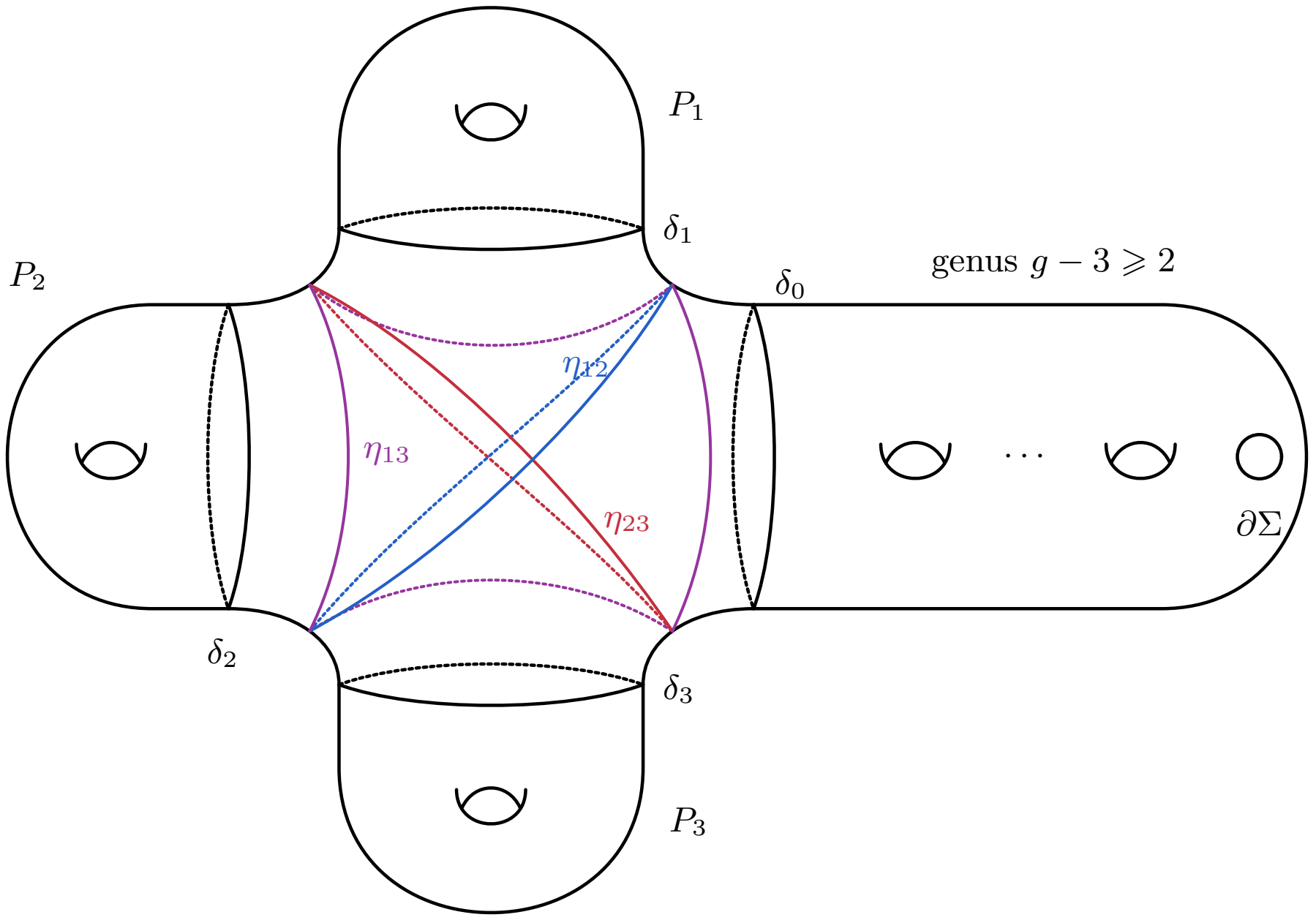

For $f\in\I$, write
$\alpha(f)=v_1+v_{\mathrm{out}}$, where
$v_1\in H_1(P_1;\Z)$ and $v_{\mathrm{out}}$ belongs to
the homology summand of its complement.
Lemma~\ref{lem:local-alpha}, applied to a one-boundary
subsurface of that complement with the same genus,
gives $u$ supported outside $P_1$ with
$\alpha(u)=v_{\mathrm{out}}$.
Since $u$ commutes with $T_{\delta_1}$, the element
$f_0:=fu^{-1}$ satisfies
\[
 \alpha(f_0)=v_1,\qquad
 [f,T_{\delta_1}]=[f_0,T_{\delta_1}].
\]
The genus-two subsurface bounded by $\eta_{13}$ contains
$P_1$ and is disjoint from $\delta_2$.
Applying the lemma there gives $u_2$ with
$\alpha(u_2)=v_1$ that commutes with $T_{\delta_2}$.
Equation~\eqref{eq:correction} gives
$[f_0,T_{\delta_2}]\in\Gamma'$.
Applying the lemma in the subsurface bounded by
$\eta_{12}$ likewise gives
$[f_0,T_{\delta_3}]\in\Gamma'$.
Consequently, $f_0\Gamma'$ commutes with every factor
of \eqref{eq:lantern} except possibly $T_{\delta_1}\Gamma'$.
Conjugating the relation by $f_0\Gamma'$ and cancelling
the other boundary twists gives
$[f_0,T_{\delta_1}]\in\Gamma'$, as required.

\emph{The side away from $\partial\Sigma$ has genus $g-1$.}
These are the only remaining curves.
Inside that side choose subsurfaces of genera $1,1,g-3$,
each with one boundary component, joined by a four-holed
sphere whose fourth boundary is $c$.
The lantern relation expresses $T_c$ as a product of
twists and inverse twists bounding genera
$1,1,g-3,2,g-2,g-2$.
The genus-one twists were treated above. All other curves
have genus at least $2$ on both sides because $g\geqslant5$,
so their twists were treated in the first case.
Their images, and hence the image of $T_c$, are central
in $N$. This proves \eqref{eq:centrality}.

Finally, Faes--Massuyeau--Sato prove
\begin{equation}\label{eq:fms-third}
 \gamma_3\I=[\I,K]
\end{equation}
\cite[Theorem~A]{FMS}. Hence $\gamma_3\I\leqslant\Gamma'$
and $N$ has nilpotency class at most $2$.
\end{proof}

\begin{remark}\label{rem:parity}
The proof uses only $K\leqslant\ker\alpha\leqslant\Gamma$.
In particular, it applies to every $n$, including even $n$.
The genus-one case is handled by the lantern relation and
supported corrections on the other two tori.
\end{remark}

For $x\in\I$, we have $\alpha(x^n)=n\alpha(x)$ and
therefore $x^n\in\Gamma$.
Since $N$ has nilpotency class at most $2$, its commutators
satisfy
\[
 [x\Gamma',y\Gamma']^{n^2}
 =[x^n\Gamma',y^n\Gamma']=1.
\]
These commutators generate $\I'/\Gamma'$, so
$n^2(\I'/\Gamma')=0$. We next determine this quotient.

\section{The commutator quotient}\label{sec:commutator}

Assume $g\geqslant5$. By Proposition~\ref{prop:third-term},
\begin{equation}\label{eq:degree-two-quotient}
 \I'/\Gamma'
 \cong\frac{\I'/\gamma_3\I}{\Gamma'/\gamma_3\I}.
\end{equation}
We compute this quotient in the integral image of the
second Johnson and Casson--Morita homomorphisms.

\subsection{The degree-two Torelli quotient}\label{sec:lattice}

Recall $\lambda=d_{\mathrm M}/8:K\to\Z$ from
\eqref{eq:lambda-definition}. Put
\[
 \mathcal T_{\Z}:=\Sym^2(\Lambda^2H)/j(\Lambda^4H),
 \qquad\mathcal T_{\Q}:=\mathcal T_{\Z}\otimes\Q,
\]
where
\begin{equation}\label{eq:plucker}
 j(x\wedge y\wedge z\wedge t)
 =(x\wedge y)(z\wedge t)-(x\wedge z)(y\wedge t)
 +(x\wedge t)(y\wedge z).
\end{equation}
The group $\mathcal T_{\Z}$ is free: for each set of four
distinct basis labels, its single relation eliminates
one of the three corresponding monomials with coefficient
$1$, and different sets involve disjoint monomials.
We identify it with its image in $\mathcal T_{\Q}$.

Let $\tau_2:K\to\mathcal T_{\Q}$ be the second Johnson
homomorphism. Faes--Massuyeau--Sato
\cite[Theorem~4.2(i)]{FMS} prove that
$(\tau_2,8d'')$ is injective on $\I'/\gamma_3\I$.
Their equation~(3.5) is $d_{\mathrm M}=2d'+48d''$,
where $d'$ is determined by $\tau_2$.
Thus $(\tau_2,\lambda)$ has the same kernel there.
Set
\begin{equation}\label{eq:fms-input}
 \Phi:\I'/\gamma_3\I\xrightarrow{\cong}M,
 \qquad\Phi(f\gamma_3\I)=(\tau_2(f),\lambda(f)),
\end{equation}
where $M$ is its image in $\mathcal T_{\Q}\oplus\Z$.

If two elements of $\I$ have the same Johnson image,
they differ by an element of $K$. Since
$[\I,K]=\gamma_3\I$, replacing either entry of a commutator
by such an element does not change its class modulo
$\gamma_3\I$. The commutator is alternating and bilinear
in this quotient, so it defines a homomorphism
\[
 B:\Lambda^2V\twoheadrightarrow M,\qquad
 B(\tau(f),\tau(h))=\Phi([f,h]\gamma_3\I).
\]
It is surjective because commutators generate $\I'$.
Denote the coordinate projections on $M$ by
$\operatorname{pr}_1:M\to\mathcal T_{\Q}$ and
$\lambda_M:M\to\Z$.
For decomposable trivectors, the tree coordinate is
\cite[Proposition~3.2]{FMS}
\begin{equation}\label{eq:bracket}
 \begin{aligned}
 &\operatorname{pr}_1B(x_1\wedge x_2\wedge x_3,y_1\wedge y_2\wedge y_3)\\
 &\qquad=\sum_{i,j\in\Z/3\Z}\omega(x_i,y_j)
 (x_{i+1}\wedge x_{i+2})(y_{j+1}\wedge y_{j+2}),
 \end{aligned}
\end{equation}
with cyclic subscripts. In particular,
$\tau_2(\I')\subseteq\mathcal T_{\Z}$.
Morita's formula~\eqref{eq:morita-commutator} becomes
\begin{equation}\label{eq:contraction-bracket}
 \lambda_M(B(v,w))=\omega(Cv,Cw).
\end{equation}
It also shows that $\lambda_M$ is surjective.
Both coordinate formulas are equivariant under $\Sp(H)$,
with the trivial action on $\Z$. Hence $B$ is equivariant
and $B(\Lambda^2W)$ is invariant under integral symplectic
changes of basis.

Under \eqref{eq:fms-input}, the subgroup
$\Gamma'/\gamma_3\I$ has image
\begin{equation}\label{eq:spin-image}
 S:=B(\Lambda^2L).
\end{equation}
Indeed, the Johnson images of elements of $\Gamma$
range over $L$, and their commutators generate $\Gamma'$.
Thus $S$ is the image of the homomorphism
$\Lambda^2L\to M$, and is a $\Z$-submodule of $M$.
Since $\gamma_3\I\leqslant\Gamma'$, we obtain
\begin{equation}\label{eq:coordinate-quotient}
 \I'/\Gamma'\cong M/S.
\end{equation}

\subsection{The kernel of the Casson--Morita coordinate}
\begin{proposition}\label{prop:kernel}\label{prop:kernel-lattice}
For $g\geqslant5$,
\begin{equation}\label{eq:kernel}
 B(\Lambda^2W)=B(W,V)=\ker(\lambda_M:M\to\Z).
\end{equation}
Here $B(W,V)$ denotes the subgroup generated by the indicated
brackets.
\end{proposition}
\begin{proof}
Put $S_0:=B(\Lambda^2W)$ and $\overline M:=M/S_0$.
Equation~\eqref{eq:contraction-bracket} gives
$S_0\subseteq\ker\lambda_M$.
Write $\Omega_i=a_i\wedge b_i$; products of bivectors
below are taken in $\mathcal T_{\Q}$.

\emph{An integral identity.}
Let $y,z$ be basis vectors on different handles, and
choose three distinct handles $j,k,l$ outside those two.
For $h\in\{j,k,l\}$, put
\[
 P_h=(y\wedge a_h)(z\wedge b_h),\qquad
 Q_h=(y\wedge b_h)(z\wedge a_h).
\]
For distinct $r,s\in\{j,k,l\}$, the bracket formulas give
\begin{align*}
 B(y\wedge a_r\wedge a_s,z\wedge b_r\wedge b_s)&=(P_r+P_s,0),\\
 B(y\wedge a_r\wedge b_s,z\wedge b_r\wedge a_s)&=(Q_s-P_r,0).
\end{align*}
Each entry in these brackets belongs to $W$.
The relation \eqref{eq:plucker} gives
$(y\wedge z)\Omega_j=P_j-Q_j$. Therefore
\begin{equation}\label{eq:auxiliary}
 \begin{aligned}
 \bigl((y\wedge z)\Omega_j,0\bigr)
 ={}&B(y\wedge a_j\wedge a_l,z\wedge b_j\wedge b_l)\\
 &-B(y\wedge a_k\wedge a_l,z\wedge b_k\wedge b_l)\\
 &-B(y\wedge a_k\wedge b_j,z\wedge b_k\wedge a_j)\ \in S_0.
 \end{aligned}
\end{equation}
Indeed, the first coordinate on the right is
$(P_j+P_l)-(P_k+P_l)-(Q_j-P_k)=P_j-Q_j$,
and its second coordinate is zero.

\emph{Brackets with one entry in $W$.}
First take an elementary trivector $t$ whose three basis
vectors lie on different handles.
For an elementary trivector $w$ with $Cw=0$, both entries
belong to $W$, so $B(t,w)\in S_0$.
Every other elementary trivector can be written, up to
sign, as $w=\Omega_i\wedge x$, where $x$ is a basis
vector outside handle $i$.
Choose a handle $j$ outside the handles of $t$ and $x$;
at most four handles are excluded.
Since $Cw=C(\Omega_j\wedge x)=x$, their difference
belongs to $W$, and
\[
 B(t,w)\equiv B(t,\Omega_j\wedge x)\pmod{S_0}.
\]
The only possible nonzero mixed pairing on the right
is between $x$ and one of the three entries of $t$.
Thus the bracket is either zero or
$\pm((y\wedge z)\Omega_j,0)$, where $y,z$ are the
other two entries. Equation~\eqref{eq:auxiliary} puts
this bracket in $S_0$. Since elementary trivectors
generate $V$, we have $B(t,V)\subseteq S_0$.

To pass from these $t$ to all of $W$, group the remaining
elementary trivectors as $\Omega_i\wedge x$, according
to their contraction $x$.
Within each such group, the kernel of $C$ is generated
by the differences $(\Omega_i-\Omega_j)\wedge x$.
Together with the elementary trivectors considered above,
these differences generate $W$ over $\Z$.
They satisfy the integral identity
\begin{equation}\label{eq:primitive}
 \begin{aligned}
 (\Omega_i-\Omega_j)\wedge x
 ={}&(a_i+a_j)\wedge(b_i-b_j)\wedge x\\
 &+a_i\wedge b_j\wedge x-a_j\wedge b_i\wedge x.
 \end{aligned}
\end{equation}
The pairs $(a_i+a_j,b_i)$ and $(b_i-b_j,a_j)$ form a
symplectic basis on handles $i,j$ and are orthogonal
to all other handles. Relative to this basis, the first
term on the right has its three entries on different
handles. The last two terms have this property in the
original basis. Since $B$ is equivariant and $S_0$ is
invariant under integral symplectic changes of basis,
each term has its bracket with $V$ in $S_0$. Hence
\begin{equation}\label{eq:mixed}
 B(W,V)=S_0.
\end{equation}

\emph{The quotient by $S_0$.}
Changing either entry of $B$ by an element of $W$
does not change its image in $\overline M$.
Since $C:V\to H$ is surjective, $B$ therefore induces
a surjective homomorphism
\[
 \overline B:\Lambda^2H\twoheadrightarrow\overline M,
 \qquad\overline B(Cv,Cw)=B(v,w)+S_0.
\]
Let $x,y$ be basis vectors on different handles.
Choose distinct handles $p,q$ outside those two.
The lifts $\Omega_p\wedge x$ and $\Omega_q\wedge y$
have zero bracket, so $\overline B(x,y)=0$.
For $i\ne j$, choose distinct $p,q$ outside handles
$i,j$. The lifts
\[
 \Omega_p\wedge(a_i+a_j),\qquad
 \Omega_q\wedge(b_i-b_j)
\]
also have zero bracket, since their only possible mixed
pairing is $\omega(a_i+a_j,b_i-b_j)=0$.
Expanding and using the vanishing just proved gives
\[
 \overline B(a_i,b_i)=\overline B(a_j,b_j).
\]
Thus $\overline M$ is generated by
$z:=\overline B(a_1,b_1)$.
The homomorphism induced by $\lambda_M$ sends $z$ to
$1$ by \eqref{eq:contraction-bracket}.
Every element of $\overline M$ is an integer multiple
of $z$, and such a multiple has zero image only when
the integer is zero. Hence $\overline M\to\Z$ is an
isomorphism, proving $S_0=\ker\lambda_M$.
Together with \eqref{eq:mixed}, this proves
\eqref{eq:kernel} over $\Z$.
The choices of handles above are possible for $g\geqslant5$.
\end{proof}

\subsection{The Casson--Morita isomorphism}

\begin{lemma}\label{lem:spin-lattice}
For $g\geqslant5$, the subgroup in \eqref{eq:spin-image}
satisfies
\begin{equation}\label{eq:spin-lattice}
 S=\lambda_M^{-1}(n^2\Z).
\end{equation}
\end{lemma}
\begin{proof}
For $v,w\in L$, we have $Cv,Cw\in nH$.
Equation~\eqref{eq:contraction-bracket} gives
$\lambda_M(B(v,w))\in n^2\Z$, so
$S\subseteq\lambda_M^{-1}(n^2\Z)$.
Conversely, $W\subseteq L$ and
Proposition~\ref{prop:kernel-lattice} give
$\ker\lambda_M\subseteq S$.
Choose
\[
 p=a_1\wedge a_2\wedge b_2,\qquad
 q=b_1\wedge a_3\wedge b_3.
\]
Then $Cp=a_1$ and $Cq=b_1$, so
$s_0:=B(np,nq)\in S$ satisfies $\lambda_M(s_0)=n^2$.
If $z\in\lambda_M^{-1}(n^2\Z)$, write
$\lambda_M(z)=an^2$ with $a\in\Z$. We have
$z-as_0\in\ker\lambda_M\subseteq S$, and therefore
$z\in S$. This proves the reverse inclusion.
\end{proof}

\begin{proposition}\label{prop:full-boundary-derived}
Let $g\geqslant5$. The homomorphism
\begin{equation}\label{eq:full-boundary-derived}
 d:\I'/\Gamma'\longrightarrow\Z/n^2\Z,
 \qquad d(f\Gamma')=\lambda(f)\pmod{n^2},
\end{equation}
is an isomorphism.
\end{proposition}
\begin{proof}
By \eqref{eq:spin-image} and Lemma~\ref{lem:spin-lattice},
$\lambda(\Gamma')\subseteq n^2\Z$, so $d$ is well defined.
If $d(f\Gamma')=0$, then
\[
 \Phi(f\gamma_3\I)\in\lambda_M^{-1}(n^2\Z)
 =S=\Phi(\Gamma'/\gamma_3\I).
\]
The injectivity of $\Phi$ gives
$f\in\Gamma'\gamma_3\I=\Gamma'$.
Thus $d$ is injective. Since $\lambda(\I')=\Z$,
its reduction modulo $n^2$ is surjective, and so is $d$.
\end{proof}

\begin{corollary}\label{cor:boundary-commutator}\label{cor:boundary-derived}
Choose $f,h\in\I$ with
$\tau(f)=a_1\wedge a_2\wedge b_2$ and
$\tau(h)=b_1\wedge a_3\wedge b_3$. Then
$[f,h]\Gamma'$ generates $\I'/\Gamma'$ and has order
$n^2$. Moreover,
\[
 \Gamma'=\{x\in\I':d_{\mathrm M}(x)\in8n^2\Z\}.
\]
\end{corollary}
\begin{proof}
Equation~\eqref{eq:morita-commutator} gives
$\lambda([f,h])=\omega(a_1,b_1)=1$.
Proposition~\ref{prop:full-boundary-derived} therefore
sends $[f,h]\Gamma'$ to the generator $1\bmod n^2$.
The same proposition characterizes membership in
$\Gamma'$ by $\lambda(x)\in n^2\Z$;
using $d_{\mathrm M}=8\lambda$ gives the last assertion.
\end{proof}

\section{Splitting the torsion sequence}\label{sec:splitting}

We determine the extension~\eqref{eq:torsion-extension}
using $\beta$ and the integral homomorphism
$\lambda=d_{\mathrm M}/8:K\to\Z$.

\begin{proposition}\label{prop:split-torsion}
For $g\geqslant5$, the map
\begin{equation}\label{eq:split-torsion}
 \begin{split}
 \Psi:K/\Gamma'&\longrightarrow
 B_{\leqslant2}\oplus\Z/n^2\Z,\quad(\cong K/\I'\oplus \I'/\Gamma')\\
 k\Gamma'&\longmapsto
 (\beta(k),\lambda(k)\bmod n^2)
 \end{split}
\end{equation}
is an isomorphism. In particular the sequence \eqref{eq:torsion-extension} splits.
\end{proposition}
\begin{proof}
The map is well defined since $\Gamma'\leqslant\I'$
and $\lambda(\Gamma')\subseteq n^2\Z$.
If $\Psi(k\Gamma')=0$, then $\beta(k)=0$ gives
$k\in\I'$ by~\eqref{eq:johnson-torsion}.
Since $\lambda(k)\equiv0\pmod{n^2}$,
Proposition~\ref{prop:full-boundary-derived}
gives $k\in\Gamma'$. Thus $\Psi$ is injective.

For surjectivity, choose $c\in\I'$ with $\lambda(c)=1$.
Then $\beta(c)=0$. Given $(b,t)$ in the target, choose
$k\in K$ with $\beta(k)=b$ and an integer lift
$\widetilde t$ of $t$. We have
\[
 \Psi\bigl(kc^{\widetilde t-\lambda(k)}\Gamma'\bigr)=(b,t).
\]
This proves surjectivity for every $n$.
\end{proof}

\begin{proof}[Proof of the bordered assertions of
Theorems~\ref{thm:main} and~\ref{thm:derived}]
Propositions~\ref{prop:third-term} and
\ref{prop:full-boundary-derived} give the assertions
about commutator subgroups.
Since $L$ is free abelian, the exact sequence
\eqref{eq:free-sequence} splits.
Proposition~\ref{prop:split-torsion} then gives
\[
 \Gamma^{\ab}\cong L\oplus B_{\leqslant2}\oplus\Z/n^2\Z.
\]
The rank of $L$ is $\binom{2g}{3}$ by
Corollary~\ref{cor:free-quotient}, and the dimension of
$B_{\leqslant2}$ is $2g^2+g+1$ by
\eqref{eq:boolean-dimensions}. This proves
\eqref{eq:boundary-answer}.
\end{proof}

\begin{proof}[Proof of Corollary~\ref{cor:intro-cyclic-generator}]
The bounding-pair formula~\eqref{eq:bp-johnson} gives
$C\tau(f)=\pm a$ and $C\tau(h)=\pm b$.
Thus \eqref{eq:morita-commutator} gives
$\lambda([f,h])=\pm\omega(a,b)=\pm1$.
Also $[f,h]\in\I'\leqslant K\leqslant\Gamma$ and
$\beta([f,h])=0$.
Proposition~\ref{prop:split-torsion} sends its class to
$(0,\pm1)$ in $B_{\leqslant2}\oplus\Z/n^2\Z$.
Since $K/\Gamma'$ is a direct summand of $\Gamma^{\ab}$
by \eqref{eq:free-sequence}, this class has order $n^2$
and generates a direct summand.
\end{proof}

Representatives for the other summands can be chosen
as follows. Choose a basis of $L$ and lift it to elements
of $\Gamma$; their classes give a free summand of
$\Gamma^{\ab}$.
Choose $c\in\I'$ with $\lambda(c)=1$.
For a basis $e_1,\ldots,e_{2g^2+g+1}$ of $B_{\leqslant2}$,
choose $k_j\in K$ with $\beta(k_j)=e_j$.
Then
\[
 \Psi\bigl(k_jc^{-\lambda(k_j)}\Gamma'\bigr)=(e_j,0).
\]
These classes generate the $B_{\leqslant2}$ summand.
Together with the free classes and the commutator in
Corollary~\ref{cor:intro-cyclic-generator}, they generate
$\Gamma^{\ab}$.

\section{Capping and the closed surface}\label{sec:capping}

Assume $g\geqslant5$ and $n\mid g-1$. We compute
$\I_g'/\Gamma_g(r)'$ and $\K_g/\Gamma_g(r)'$ by capping.
The Birman exact sequence is
\begin{equation}\label{eq:birman}
 1\longrightarrow\mathcal P\longrightarrow\I_{g,1}
 \xrightarrow{\operatorname{cap}}\I_g\longrightarrow1,
\end{equation}
where $\mathcal P\cong\pi_1(UT\Sigma_g)$ is the
disk-pushing subgroup. Choose the fiber generator
$t=T_\partial$ to be the positive boundary twist.
The Johnson homomorphism on $\mathcal P$ factors as
\begin{equation}\label{eq:push-johnson}
 \mathcal P\longrightarrow H\xrightarrow{\iota}V,
\end{equation}
with the first map induced by projection to $\Sigma_g$
and abelianization~\cite[Lemma~2.8]{Putman}.
For $p\in\mathcal P$, equations
\eqref{eq:push-johnson} and \eqref{eq:contraction-iota}
give $C\tau_{g,1}(p)\in(g-1)H\subseteq nH$.
Thus $\mathcal P\leqslant\Gamma_{g,1}(r)$.
If $f\in\Gamma_g(r)$ and $\widetilde f\in\I_{g,1}$
is a lift, then $C_n\tau_{g,1}(\widetilde f)=C_n\tau_g(f)=0$.
Hence $\widetilde f\in\Gamma_{g,1}(r)$, proving that
capping restricts to an exact sequence
\begin{equation}\label{eq:spin-birman}
 1\longrightarrow\mathcal P\longrightarrow\Gamma_{g,1}(r)
 \xrightarrow{\operatorname{cap}}\Gamma_g(r)\longrightarrow1.
\end{equation}

\begin{lemma}\label{lem:push-derived}
The disk-pushing subgroup satisfies
\[
 \mathcal P\cap\K_{g,1}=\mathcal P'\langle t\rangle,
 \qquad
 \mathcal P\cap\I_{g,1}'=\mathcal P'\langle t^2\rangle.
\]
Moreover,
\begin{equation}\label{eq:boundary-invariants}
 \beta(t)=\Arf,\qquad
 \lambda(t)=\frac{g(g-1)}2,\qquad
 \lambda(t^2)=g(g-1).
\end{equation}
\end{lemma}
\begin{proof}
Quotienting $\mathcal P$ by $\langle t\rangle$ gives
$\pi_1(\Sigma_g)$. The kernel of its abelianization is
the commutator subgroup, whose inverse image in
$\mathcal P$ is $\mathcal P'\langle t\rangle$.
Thus
\[
 \ker(\mathcal P\to H)=\mathcal P'\langle t\rangle.
\]
Since $\iota$ is injective, \eqref{eq:push-johnson}
gives $\mathcal P\cap\K_{g,1}=\mathcal P'\langle t\rangle$.
Now $\mathcal P'\leqslant\I_{g,1}'$ and
$\beta(t)=\Arf$ by \eqref{eq:bcj-twist}.
For $p'\in\mathcal P'$ and $j\in\Z$, we have
\[
 \beta(p't^j)=j\Arf.
\]
As $\Arf$ is nonzero of order $2$, this vanishes
exactly when $j$ is even. Equation~\eqref{eq:johnson-torsion}
therefore gives
$\mathcal P\cap\I_{g,1}'=\mathcal P'\langle t^2\rangle$.
Finally, $t$ bounds genus $g$, so the normalization
$d_{\mathrm M}(t)=4g(g-1)$ gives the two values of
$\lambda$ in \eqref{eq:boundary-invariants}.
\end{proof}

\begin{proposition}\label{prop:closed-derived}
Put $m:=\gcd(n^2,g-1)$. Then
\[
 \I_g'/\Gamma_g(r)'\cong\Z/m\Z,
 \qquad\gamma_3\I_g\leqslant\Gamma_g(r)'.
\]
\end{proposition}
\begin{proof}
The surjections in \eqref{eq:birman} and
\eqref{eq:spin-birman} map commutator subgroups onto
commutator subgroups. Hence
\[
 \I_g'/\Gamma_g(r)'
 \cong\frac{\I_{g,1}'}
 {\Gamma_{g,1}(r)'(\mathcal P\cap\I_{g,1}')}.
\]
Since $\mathcal P\leqslant\Gamma_{g,1}(r)$, we have
$\mathcal P'\leqslant\Gamma_{g,1}(r)'$.
Lemma~\ref{lem:push-derived} therefore shows that capping
kills precisely the subgroup generated by the class of
$t^2$ in $\I_{g,1}'/\Gamma_{g,1}(r)'$.
Under Proposition~\ref{prop:full-boundary-derived}, this
class has value $\lambda(t^2)=g(g-1)\bmod n^2$. Thus
\[
 \I_g'/\Gamma_g(r)'
 \cong(\Z/n^2\Z)/\langle g(g-1)+n^2\Z\rangle
 \cong\Z/\gcd(n^2,g(g-1))\Z.
\]
As $n\mid g-1$, we have $\gcd(n,g)=1$, and consequently
$\gcd(n^2,g(g-1))=\gcd(n^2,g-1)=m$.
This proves the isomorphism.
Finally, capping maps $\gamma_3\I_{g,1}$ onto
$\gamma_3\I_g$, so Proposition~\ref{prop:third-term}
gives $\gamma_3\I_g\leqslant\Gamma_g(r)'$.
\end{proof}

To compute the whole torsion subgroup, the class of $t$
must be retained together with both coordinates in
\eqref{eq:split-torsion}.

\begin{proposition}\label{prop:closed-torsion}
For $m=\gcd(n^2,g-1)$,
\[
 \K_g/\Gamma_g(r)'
 \cong B_{\leqslant2}^0\oplus\Z/m\Z
\]
as an abstract abelian group.
\end{proposition}
\begin{proof}
We first show that capping maps $\K_{g,1}$ onto $\K_g$.
Given $f\in\K_g$, choose a lift $\widetilde f\in\I_{g,1}$.
Since $\tau_g(f)=0$, its Johnson image belongs to
$\iota(H)$. Equation~\eqref{eq:push-johnson} gives
$p\in\mathcal P$ with
$\tau_{g,1}(p)=\tau_{g,1}(\widetilde f)$.
Then $\widetilde f p^{-1}\in\K_{g,1}$ is another lift of $f$.
Together with \eqref{eq:spin-birman}, this gives
\[
 \K_g/\Gamma_g(r)'
 \cong\frac{\K_{g,1}}
 {\Gamma_{g,1}(r)'(\mathcal P\cap\K_{g,1})}.
\]
By Lemma~\ref{lem:push-derived}, the last intersection
is $\mathcal P'\langle t\rangle$.
Since $\mathcal P'\leqslant\Gamma_{g,1}(r)'$, this quotient
is obtained from $\K_{g,1}/\Gamma_{g,1}(r)'$ by killing
the class of $t$.
Put $c_g:=g(g-1)/2$.
Proposition~\ref{prop:split-torsion} and
\eqref{eq:boundary-invariants} identify the quotient with
\begin{equation}\label{eq:closed-joint-relation}
 \frac{B_{\leqslant2}\oplus\Z/n^2\Z}
 {\langle(\Arf,c_g\bmod n^2)\rangle}.
\end{equation}
Choose an $\F_2$-linear complement $B'$ to
$\langle\Arf\rangle$ in $B_{\leqslant2}$.
Restriction identifies $B'$ with $B_{\leqslant2}^0$,
and the relation in \eqref{eq:closed-joint-relation}
lies entirely in the other factor.
That factor has additive presentation
\[
 \langle a,z\mid 2a=0,\ n^2z=0,\ a+c_gz=0\rangle.
\]
Substituting $a=-c_gz$ gives
\[
 \langle z\mid n^2z=0,\ g(g-1)z=0\rangle
 \cong\Z/\gcd(n^2,g(g-1))\Z
 =\Z/m\Z.
\]
Therefore the quotient is
$B_{\leqslant2}^0\oplus\Z/m\Z$, for either parity of $n$.
\end{proof}

\begin{proof}[Proof of the closed assertions of
Theorems~\ref{thm:main} and~\ref{thm:derived}]
Proposition~\ref{prop:closed-derived} gives the assertions
about commutator subgroups.
The Johnson homomorphism gives an exact sequence
\[
 0\longrightarrow\K_g/\Gamma_g(r)'
 \longrightarrow\Gamma_g(r)^{\ab}
 \xrightarrow{\tau_g}\tau_g(\Gamma_g(r))\longrightarrow0.
\]
By Corollary~\ref{cor:free-quotient}, the last group is
free of rank $\binom{2g}{3}-2g$, so the sequence splits.
Proposition~\ref{prop:closed-torsion} identifies its kernel
with $B_{\leqslant2}^0\oplus\Z/m\Z$.
Using \eqref{eq:boolean-dimensions} gives
\eqref{eq:closed-answer}.
\end{proof}

When $n$ is even, $\lambda(t)$ need not vanish modulo $m$.
Thus $\lambda\bmod m$ is not asserted to descend from
$\K_{g,1}$ to $\K_g$. The commutator quotient is detected
using lifts in $\I_{g,1}'$, where $t^2$ is the relation;
the full torsion subgroup is computed by the joint
relation~\eqref{eq:closed-joint-relation}.

The closed Chillingworth subgroup is defined as the
capping image of $\ker(C\tau_{g,1})$
\cite[Definition~3.2]{Kosuge}. This image is
$\Gamma_g(2g-2)$. Indeed, a lift $f\in\I_{g,1}$ of
an element of $\Gamma_g(2g-2)$ has
$C\tau_{g,1}(f)=(g-1)x$ for some $x\in H$.
By \eqref{eq:push-johnson}, there is $p\in\mathcal P$
with $\tau_{g,1}(p)=\iota(x)$. Then
$fp^{-1}\in\ker(C\tau_{g,1})$ is another lift.
The reverse inclusion follows from
Proposition~\ref{prop:spin-kernel}. When $g-1$ is odd,
\eqref{eq:effective-level} also identifies this group
with $\Gamma_g(g-1)$.

\section{Genera three and four}\label{sec:low-genus}
For $g=1$, $\Gamma_{1,1}(r)=\I_{1,1}\cong\Z$,
generated by the boundary twist.
For $g=2$, capping maps $\Gamma_{2,1}(r)$ onto $\I_2$,
so $H_1(\Gamma_{2,1}(r);\Z)$ has infinite rank
by Mess's theorem~\cite{Mess}.

For surfaces with one boundary component, we show that
$\I_{g,1}'/\Gamma_{g,1}(r)'$ is not cyclic when
$g=3$ and $n>1$, or when $g=4$ and $3\mid n$.
We use the degree-two description of Faes--Massuyeau--Sato
and the commutator formulas in Section~\ref{sec:lattice},
which hold for every $g\geqslant3$.
These arguments do not assume
$\gamma_3\I_{g,1}\leqslant\Gamma_{g,1}(r)'$.

\subsection{Genus three}

\begin{proposition}\label{prop:genus-three-obstruction}
Let $g=3$, $n>1$, and let $p\mid n$ be prime. Then
\[
 \dim_{\F_p}
 \left((\I_{3,1}'/\Gamma_{3,1}(r)')\otimes\F_p\right)
 \geqslant15.
\]
In particular, $\I_{3,1}'/\Gamma_{3,1}(r)'$ is not
isomorphic to $\Z/n^2\Z$.
\end{proposition}
\begin{proof}
Put $E:=\I_{3,1}'/\gamma_3\I_{3,1}$.
By \cite[Theorem~4.2(i)]{FMS}, this group is free of rank
\[
 \dim\Sym^2(\Lambda^2\Q^6)-\dim\Lambda^4\Q^6+1
 =\frac{15\cdot16}{2}-15+1=106.
\]
Hence $E\otimes\F_p$ has dimension $106$.
The reduction of $C:V\to H$ modulo $p$ is surjective;
write $W_p$ for its kernel. Then
\[
 \dim_{\F_p}W_p=\binom63-6=14.
\]
For $f\in\Gamma_{3,1}(r)$, we have
$C\tau(f)\in nH\subseteq pH$.
Thus its Johnson image modulo $p$ belongs to $W_p$.
Since $\Gamma_{3,1}(r)'$ is generated by commutators
of such elements, its image in $E\otimes\F_p$ is
contained in the image of
\[
 \Lambda^2W_p\longrightarrow E\otimes\F_p
\]
under the commutator map. This image has dimension at
most $\binom{14}{2}=91$.
It follows that
\[
 \dim_{\F_p}\left(
 \frac{\I_{3,1}'}
 {\Gamma_{3,1}(r)'\gamma_3\I_{3,1}}
 \otimes\F_p\right)\geqslant106-91=15.
\]
This vector space is a quotient of
$(\I_{3,1}'/\Gamma_{3,1}(r)')\otimes\F_p$,
which proves the bound.
Since $(\Z/n^2\Z)\otimes\F_p$ has dimension $1$,
the two abelian groups cannot be isomorphic.
\end{proof}

\subsection{Genus four}

In genus $4$, reduction modulo $3$ detects additional
commutator classes. We use the integral tensor group
$\mathcal T_{\Z}$ from Section~\ref{sec:lattice}.
For the next proof write $H_{\F_3}=H\otimes\F_3$ and
$H_i=\F_3a_i\oplus\F_3b_i$ for $1\leqslant i\leqslant4$.

\begin{proposition}\label{prop:genus-four-obstruction}
Let $g=4$ and $3\mid n$. Then
\[
 \dim_{\F_3}
 \left((\I_{4,1}'/\Gamma_{4,1}(r)')\otimes\F_3\right)
 \geqslant17.
\]
In particular, $\I_{4,1}'/\Gamma_{4,1}(r)'$ is not
isomorphic to $\Z/n^2\Z$.
\end{proposition}
\begin{proof}
Wedging two bivectors defines a map
\[
 \operatorname{alt}:\mathcal T_{\Z}\otimes\F_3
 \longrightarrow\Lambda^4H_{\F_3},
 \qquad (x\wedge y)(z\wedge t)\longmapsto x\wedge y\wedge z\wedge t.
\]
It is well defined because the image of
$j(x\wedge y\wedge z\wedge t)$ under this rule is
$3x\wedge y\wedge z\wedge t$, which vanishes modulo $3$.
Let $P$ project $\Lambda^4H_{\F_3}$ onto the summand
with one entry on each of the four handles, identified with
\[
 Q_4:=H_1\otimes H_2\otimes H_3\otimes H_4,
 \qquad\dim_{\F_3}Q_4=16.
\]
Since $\tau_2(\I_{4,1}')\subseteq\mathcal T_{\Z}$,
reduction modulo $3$ gives a homomorphism
$\theta:=P\operatorname{alt}\tau_2:\I_{4,1}'\to Q_4$.

Define
\[
 C_6(x_1\wedge\cdots\wedge x_6)
 :=\sum_{i<j}(-1)^{i+j-1}\omega(x_i,x_j)
 x_1\wedge\cdots\widehat{x_i}\cdots\widehat{x_j}\cdots\wedge x_6.
\]
For decomposable $v,w\in\Lambda^3H_{\F_3}$, divide the
terms in $C_6(v\wedge w)$ according to the paired entries.
Pairs within $v$ give $(Cv)\wedge w$, pairs within $w$
give $v\wedge(Cw)$, and mixed pairs give the alternating
image of \eqref{eq:bracket}. By bilinearity, for all $v,w$,
\begin{equation}\label{eq:genus-four-contraction}
 \operatorname{alt}(\operatorname{pr}_1B(v,w))
 =C_6(v\wedge w)-(Cv)\wedge w-v\wedge(Cw).
\end{equation}
Also $PC_6=0$: on a basis six-vector, a nonzero contraction
removes both basis vectors of one handle. The resulting
four-vector has no entry on that handle, and its projection
under $P$ is zero.

For $f,h\in\Gamma_{4,1}(r)$, the vectors
$C\tau(f),C\tau(h)$ vanish modulo $3$ because $3\mid n$.
Applying $P$ to \eqref{eq:genus-four-contraction} gives
$\theta([f,h])=0$.
Morita's formula also gives $\lambda([f,h])=0\bmod3$.
Thus the homomorphism
\[
 (\theta,\lambda\bmod3):\I_{4,1}'\longrightarrow Q_4\oplus\F_3
\]
vanishes on $\Gamma_{4,1}(r)'$.

To prove surjectivity, choose $x_i\in\{a_i,b_i\}$ and
$x_1^*\in\Z a_1\oplus\Z b_1$ with
$\omega(x_1,x_1^*)=1$. Put
\[
 v=x_1\wedge x_2\wedge x_1^*,\qquad
 w=x_3\wedge x_4\wedge x_1.
\]
The only nonzero mixed pairing is between $x_1^*$ and
$x_1$, so the tree coordinate of $B(v,w)$ is
$-(x_1\wedge x_2)(x_3\wedge x_4)$.
Since $Cw=0$, its $\lambda_M$-coordinate is zero.
Choose $f,h\in\I_{4,1}$ with $\tau(f)=v$ and $\tau(h)=w$.
As the $x_i$ vary, the images of these commutators are,
up to sign, the sixteen basis vectors of $Q_4\oplus0$.
Next take Johnson lifts of
\[
 p=a_1\wedge a_2\wedge b_2,\qquad
 q=b_1\wedge a_3\wedge b_3.
\]
Their commutator has $\lambda$-value $1$ and tree
coordinate $(a_2\wedge b_2)(a_3\wedge b_3)$.
Its image under $\theta$ is zero, so it maps to $(0,1)$.
These images span $Q_4\oplus\F_3$.
The map therefore induces a surjection from
$(\I_{4,1}'/\Gamma_{4,1}(r)')\otimes\F_3$
onto this $17$-dimensional vector space, proving the bound.
\end{proof}

These propositions do not determine the full
abelianizations in genera $3$ and $4$. In particular,
the genus-four cases with $3\nmid n$ are not settled by
the argument above. The case $n=1$ gives the ordinary
Torelli group and is covered by Johnson's theorem.

\section*{Acknowledgements}

The author thanks Yi Liu for suggesting this problem and for helpful discussions.

\section*{Declaration on use of AI tools}
AI-based tools were used for: (\romannumeral1) the organization and wording of the text; and (\romannumeral2) the exploration of proof strategies for Proposition 6.1 and the verification of mathematical arguments and computations.
The author independently verified and revised all AI-assisted mathematical arguments and takes full responsibility for the mathematical correctness and final content of the manuscript.

\par\bigskip
\begingroup
\small
\setlength{\parindent}{2em}

\textsc{School of Mathematical Sciences,
Beijing Normal University, Beijing, 100875, The People's Republic of China}\par
\textit{Email address:}
\href{mailto:202011130077@mail.bnu.edu.cn}
     {\texttt{202011130077@mail.bnu.edu.cn}}\par

\endgroup

\end{document}